\pdfoutput=1
\documentclass[11pt]{article}
\usepackage[T1]{fontenc}
\usepackage{lmodern}

\usepackage[margin=1in]{geometry}
\usepackage{amsmath,amssymb,amsthm,mathtools}
\usepackage{enumitem}
\usepackage{pgfplots}
\pgfplotsset{compat=1.17}
\usepackage{color}
\usepackage[colorlinks=true,linkcolor=blue,citecolor=blue,urlcolor=blue]{hyperref}
\hypersetup{pdftitle={Two-Stage Mixed-Integer Quadratic Programming in Polynomial Time},pdfauthor={Alberto Del Pia}}
\makeatletter\g@addto@macro{\UrlBreaks}{\do\-}\makeatother

\newtheorem{theorem}{Theorem}
\newtheorem{lemma}[theorem]{Lemma}
\newtheorem{proposition}[theorem]{Proposition}
\newtheorem{corollary}[theorem]{Corollary}
\theoremstyle{definition}
\newtheorem{algo}{Algorithm}

\newcommand{\R}{\mathbb{R}}
\newcommand{\Z}{\mathbb{Z}}
\newcommand{\Q}{\mathbb{Q}}
\newcommand{\cJ}{\mathcal{J}}
\newcommand{\cL}{\mathcal{L}}
\DeclareMathOperator{\sign}{sign}

\title{Two-Stage Mixed-Integer Quadratic Programming\\ in Polynomial Time}
\author{Alberto Del Pia\thanks{Department of Industrial and Systems Engineering \& Wisconsin Institute for Discovery, University of Wisconsin--Madison. Email: \texttt{delpia@wisc.edu}.}}
\date{November 10, 2026}

\begin{document}
\maketitle

\begin{abstract}
Lenstra showed that mixed-integer linear programming can be solved in polynomial time when the number of integer variables is fixed and the number of continuous variables is arbitrary. Mixed-integer quadratic programming (MIQP) can be solved in polynomial time when the total number of variables is fixed. However, since MIQP is NP-hard already without integer variables, the analog of Lenstra's result fails unless P${}={}$NP. Hence, polynomial-time algorithms for MIQP with an arbitrary number of continuous variables must exploit additional structure. We study MIQPs with two-stage structure, which underlies Benders decomposition and two-stage stochastic programming. A few linking variables, including all integer variables, are coupled with many blocks of continuous variables that interact only through them, both in the constraints and in the objective. We give a polynomial-time algorithm for MIQP with an arbitrary quadratic objective when the number of linking variables and the block size are fixed, while the number of blocks, and hence the number of variables, is arbitrary. The case with no blocks is MIQP in fixed dimension. None of our assumptions can be dropped. Unless P${}={}$NP, the number of linking variables and the block size must both be fixed, the blocks cannot interact or contain integer variables, and large blocks cannot be allowed even if their objective is linear. Unless FPT${}={}$W[1], the exponent of the running time must depend on both the number of linking variables and the block size.
\end{abstract}


\section{Introduction}\label{sec:intro}

Lenstra showed that mixed-integer linear programming can be solved in polynomial time when the number of integer variables is fixed and the number of continuous variables is arbitrary~\cite[Section~5]{Len83}. \emph{Mixed-integer quadratic programming} (MIQP), where the linear objective is replaced by a quadratic function, can be solved in polynomial time when the total number of variables is fixed~\cite{AriHil26e} (see also~\cite{Wei26} for the pure integer case). However, unless P${}={}$NP, MIQP cannot be solved in polynomial time when only the number of integer variables is fixed. Indeed, an MIQP without integer variables is a continuous quadratic program, which is NP-hard already when the quadratic form has one negative eigenvalue~\cite{ParVav91}.

For convex objectives, Lenstra's result does extend: mixed-integer convex quadratic programming can be solved in polynomial time when the number of integer variables is fixed~\cite{dP25SIOPT}. For nonconvex objectives and a fixed number of integer variables, approximation algorithms are known when the quadratic form has a fixed number of negative eigenvalues~\cite{dP26Jacobi} (see also~\cite{dP18MPB,dP23bMPA,Vav92MPB}). This raises the question of which structure allows nonconvex MIQPs with a fixed number of integer variables and arbitrarily many continuous variables to be solved exactly in polynomial time.

\paragraph{Two-stage MIQP.} We study the classical two-stage, or dual block-angular, structure, in which a few linking variables are coupled with many blocks of variables that interact only through them. It underlies Benders decomposition for mixed-integer linear programming~\cite{Ben62} and two-stage stochastic programming. Under the name of two-stage stochastic integer programming, it has been studied for integer linear programs~\cite{EisHunKleKouLevOnn19,HemSch03,Kle22} and for mixed-integer linear programs~\cite{BraKouLasOrd24}. We show that MIQPs with this structure and an arbitrary quadratic objective can be solved exactly in polynomial time when the number of linking variables and the size of the blocks are fixed. We consider problems of the form
\begin{equation}\label{eq:bmiqp}
\begin{aligned}
\min\quad & x^\top Q_0 x + c_0^\top x + \sum_{i=1}^n f_i(x, y_i) \\
\text{s.t.}\quad & A_0 x \le b_0, \\
& A_i x + D_i y_i \le b_i, && i = 1, \dots, n, \\
& x \in \Z^p \times \R^q, \quad y_i \in \R^{r_i}, && i = 1, \dots, n,
\end{aligned}
\end{equation}
where
\[
f_i(x, y_i) = y_i^\top Q_i y_i + 2 x^\top R_i y_i + c_i^\top y_i, \qquad i = 1, \dots, n.
\]
All data are rational, and the matrices $Q_0$ and $Q_1, \dots, Q_n$ are symmetric. The vector $x$ consists of the \emph{linking variables}. The first $p$ of them are integral and the remaining $q$ are the \emph{continuous linking variables}, and we write $d = p + q$. For each $i$, the vector $y_i$ is the $i$th \emph{block} of continuous variables. The problem has $m_0$ linking constraints $A_0 x \le b_0$, and block $i$ has $m_i$ constraints $A_i x + D_i y_i \le b_i$. We write $m = m_0 + m_1 + \dots + m_n$. The \emph{block size} is $r = \max_i r_i$, and $r = 0$ if $n = 0$. A \emph{feasible point} of~\eqref{eq:bmiqp} is a point $(x, y_1, \dots, y_n)$ that satisfies all constraints of~\eqref{eq:bmiqp}. We denote by $L$ the binary encoding length of an instance, and \emph{polynomial time} means time polynomial in $L$.

\paragraph{Results.} Our main result is the following.

\begin{theorem}\label{thm:main}
For all fixed $d$ and $r$, there is a polynomial-time algorithm that, given an instance of~\eqref{eq:bmiqp}, either reports that~\eqref{eq:bmiqp} has no feasible point, or reports that the objective of~\eqref{eq:bmiqp} is unbounded below on its feasible points, or returns an optimal solution.
\end{theorem}

Two special cases of Theorem~\ref{thm:main} are covered by known results. If $n = 0$, then~\eqref{eq:bmiqp} is an arbitrary MIQP in $d$ variables, and Theorem~\ref{thm:main} is the polynomial-time solvability of MIQP in fixed dimension~\cite{AriHil26e}. If the objective is convex, Theorem~\ref{thm:main} follows from~\cite[Theorem~1.3]{dP25SIOPT}, and if it is linear, already from Lenstra's algorithm~\cite[Section~5]{Len83}, \cite[Section~18.4]{Sch86}. Both Lenstra's algorithm and the algorithm of~\cite{dP25SIOPT} require only $p$ to be fixed. Hence Theorem~\ref{thm:main} is new only for objectives that are not convex.

Two-stage stochastic programming with finitely many scenarios is a natural source of instances of~\eqref{eq:bmiqp}. There, $x$ is the first-stage decision, block $i$ is the recourse decision in scenario $i$, and $n$ is the number of scenarios. Theorem~\ref{thm:main} applies when there are few first-stage decisions and each scenario has a small recourse problem with an arbitrary quadratic cost, for example a cost with economies of scale, while the number of scenarios is arbitrary. Already with $d = r = 2$ and an indefinite matrix $Q_i$ in every scenario, the problem has at least $n$ negative eigenvalues, so it is covered neither by MIQP in fixed dimension nor by mixed-integer convex quadratic programming. Moreover, the optimal value of a scenario need not be piecewise quadratic in $x$ (Section~\ref{sec:main}). Theorem~\ref{thm:main} is of interest also for $p = 0$, where~\eqref{eq:bmiqp} is a continuous quadratic program with arbitrarily many negative eigenvalues, while quadratic programming is NP-hard already with one negative eigenvalue~\cite{ParVav91}.

The next result shows that, unless P${}={}$NP, none of the assumptions of Theorem~\ref{thm:main} can be dropped, even over polytopes, and convexity cannot replace the bound on the block size.

\begin{theorem}\label{thm:tight}
Unless P${}={}$NP, Theorem~\ref{thm:main} does not hold if any of the following changes is made, even when the linear constraints of~\eqref{eq:bmiqp} define a polytope.
\begin{enumerate}[label=(\roman*),leftmargin=*,widest=vii]
\item\label{it:t-d} The number $d$ of linking variables is not fixed.
\item\label{it:t-k} The block size $r$ is not fixed.
\item\label{it:t-q} The objective may contain products of variables of different blocks.
\item\label{it:t-link} The constraints may involve variables of different blocks.
\item\label{it:t-int} The blocks may contain integer variables.
\item\label{it:t-conv} There may be a single block of unbounded size whose function $f_1$ is convex in $(x, y_1)$.
\item\label{it:t-conv2} Blocks of unbounded size are allowed when the sum of $x^\top Q_0 x + c_0^\top x$ and of their functions $f_i$ is convex.
\end{enumerate}
\end{theorem}

Each part of Theorem~\ref{thm:tight} holds already for very restricted instances (Proposition~\ref{prop:strong} in Section~\ref{sec:tight}).

Parts~\ref{it:t-conv} and~\ref{it:t-conv2} of Theorem~\ref{thm:tight} concern convex objectives. For fixed $p$, MIQP can be solved in polynomial time when the objective is convex~\cite[Theorem~1.3]{dP25SIOPT}, and when it has the form~\eqref{eq:bmiqp} with fixed $q$ and $r$ (Theorem~\ref{thm:main}). Parts~\ref{it:t-conv} and~\ref{it:t-conv2} show that these two cases cannot be combined by allowing blocks of unbounded size with a convex objective in~\eqref{eq:bmiqp}.

The algorithm of Theorem~\ref{thm:main} runs in time $L^{\gamma(d,r)}$ for some function $\gamma$. The next result shows that the exponent must depend on both $d$ and $r$.

\begin{proposition}\label{thm:xp}
Unless FPT${}={}$W[1], there is no algorithm that solves the instances of~\eqref{eq:bmiqp} over polytopes in time $\omega(d,r) \cdot L^{c}$, for a computable function $\omega$ and a constant $c$. This holds already for the instances with $n = 0$, parameterized by $d$, and for the instances with $d = 0$ and $n = 1$, parameterized by $r$.
\end{proposition}

Theorem~\ref{thm:main} also applies when the block structure is not given, since it can be found by enumeration, as we explain next. The \emph{interaction graph} of an MIQP has the continuous variables as nodes, and two of them are adjacent if they appear with nonzero coefficients in a common constraint or in a common term of the objective. It is the mixed primal graph of~\cite{EibGanKnoOrd19} restricted to the continuous variables. An MIQP can be written in the form~\eqref{eq:bmiqp} with $q$ continuous linking variables and blocks of size at most $r$ if and only if deleting some $q$ nodes of its interaction graph leaves connected components with at most $r$ nodes. The \emph{vertex integrity} of a graph~\cite{BarEntSwa87} is the minimum, over all sets $W$ of nodes, of $|W|$ plus the largest number of nodes in a connected component after deleting $W$. If the interaction graph has vertex integrity at most $v$, then the MIQP can be written in the form~\eqref{eq:bmiqp} with $q \le v$ and $r \le v$, and such a form can be found by enumeration. Hence Theorem~\ref{thm:main} gives the following result.

\begin{corollary}\label{cor:vi}
For all fixed $p$ and $v$, there is a polynomial-time algorithm that, given an MIQP with $p$ integer variables whose interaction graph has vertex integrity at most $v$, either reports that it has no feasible point, or reports that its objective is unbounded below on the feasible points, or returns an optimal solution.
\end{corollary}

\paragraph{Techniques.} We outline the proof of Theorem~\ref{thm:main}. As in Benders decomposition, we fix the linking variables $x$. Problem~\eqref{eq:bmiqp} then splits into $n$ independent quadratic programs, one for each block, and it remains to minimize over $x$ the sum of $x^\top Q_0 x + c_0^\top x$ and of their optimal values. If the objective is linear, the optimal value of a block is a convex piecewise linear function of $x$, described by linear programming duality. For quadratic objectives it is in general neither convex nor described by duality. Generalized Benders decomposition~\cite{Geo72} handles nonlinear subproblems, but it requires them to be convex (see~\cite{LiGro19} for the nonconvex case). Instead, we use the following active-set argument. For each block, a minimizer is a stationary point on the affine subspace defined by a set of at most $r$ active constraints. For each such set, the stationary points form an affine subspace that depends affinely on $x$, and the objective of the block is constant on it. Hence the optimal value of block $i$ is the minimum of polynomially many quadratic functions of $x$, each defined on a polyhedron (Proposition~\ref{prop:values}). Choosing one function per block naively gives exponentially many combinations. However, the best function of each block at a point $x$ is determined by the signs of polynomially many polynomials of degree at most two in $x$, namely the differences of these functions and the inequalities that define the polyhedra. In the space of the monomials of degree at most two in $x$, these polynomials become affine functions. Their zero sets form an arrangement of polynomially many hyperplanes in fixed dimension, which has polynomially many sign vectors. Each sign vector gives one MIQP in the $d$ linking variables, which can be solved in polynomial time~\cite{AriHil26e}, or with the algorithm of Wei for integer quadratic programming~\cite{Wei26} (Appendix~\ref{app:wei}). For $p = 0$, these are quadratic programs in fixed dimension, and the classical algorithms suffice~\cite{Vav90}, \cite[Theorem~5]{dPDeyMol17MPA}. Before this, we decide whether some block is unbounded below for some value of the linking variables that extends to a feasible point. Using a ray criterion and the fact that every block has fixed dimension, this question reduces to polynomially many quadratic programs in fixed dimension and mixed-integer linear programs with $p$ integer variables.

\paragraph{Related work.} The closest results concern two-stage stochastic integer programs, which have the structure of~\eqref{eq:bmiqp} with integral blocks and a linear objective. These programs can be solved in time polynomial in the number of blocks when the number of linking variables, the block sizes, and the coefficients are bounded~\cite{EisHunKleKouLevOnn19,HemSch03,Kle22}, and the problem is then fixed-parameter tractable~\cite{EisHunKleKouLevOnn19,Kle22}. Cslovjecsek, Koutecký, Lassota, Pilipczuk, and Polak~\cite{CslKouLasPilPol25} showed that the feasibility problem of these programs remains fixed-parameter tractable when the coefficients of the linking variables are arbitrary. They also showed that it is NP-hard for blocks of constant dimension when the coefficients in the blocks are arbitrary. For mixed-integer linear programs with this structure, Brand, Koutecký, Lassota, and Ordyniak~\cite{BraKouLasOrd24} gave an algorithm for bounded coefficients, and they showed W[1]-hardness parameterized by the dimensions of the blocks, even with coefficients in $\{-1,0,1\}$. More generally, Dvo\v{r}\'ak, Eiben, Ganian, Knop, and Ordyniak~\cite{DvoEibGanKnoOrd21} determined the complexity of integer linear programs that split into small components after the deletion of a few global variables or constraints, for several encodings of the coefficients. In their terminology, the continuous linking variables of~\eqref{eq:bmiqp} play the role of fracture backdoors. For \emph{integer quadratic programming} (IQP), Eiben, Ganian, Knop, and Ordyniak~\cite{EibGanKnoOrd19} showed that structural restrictions must account for the interactions in the objective and must be combined with bounds on the domains or the coefficients.

Theorem~\ref{thm:main} differs from these results in three respects: the blocks are continuous, the objective is an arbitrary quadratic function, and the coefficients are arbitrary rational numbers. Part~\ref{it:t-int} of Theorem~\ref{thm:tight} shows that integral blocks cannot be allowed, in line with the hardness results of~\cite{CslKouLasPilPol25,DvoEibGanKnoOrd21}. Part~\ref{it:t-q} of Theorem~\ref{thm:tight} is a continuous analog of the observation of~\cite{EibGanKnoOrd19} that the interactions in the objective must be accounted for. Finally, Proposition~\ref{thm:xp} shows that, unlike the algorithms for two-stage stochastic integer programs with bounded coefficients, the algorithm of Theorem~\ref{thm:main} cannot be made fixed-parameter tractable unless FPT${}={}$W[1]. This hardness comes from the quadratic objective, and it holds already for instances without blocks or without linking variables.

\paragraph{Organization.} Section~\ref{sec:prelim} collects tools. Section~\ref{sec:block} describes the optimal value of one block as a function of the linking variables, and Section~\ref{sec:main} proves Theorem~\ref{thm:main}. Section~\ref{sec:tight} outlines the proof of Theorem~\ref{thm:tight}. The omitted proofs are in the appendix.

\section{Preliminaries}\label{sec:prelim}

For a positive integer $t$, we write $[t] = \{1, \dots, t\}$. For a matrix or vector $M$ and a set $S$ of row indices, $M_S$ denotes the submatrix of $M$ with rows indexed by $S$.

We use the following results on the polynomial-time solvability of mixed-integer linear programming with a fixed number of integer variables and of MIQP in fixed dimension.

\begin{proposition}[{\cite[Section~5]{Len83}, \cite[Section~18.4]{Sch86}}]\label{prop:milp}
For every fixed $p$, there is a polynomial-time algorithm that, given a linear function $f$ in $d \ge p$ variables and a polyhedron $P = \{x \in \R^d : Ax \le b\}$, all with rational data, either reports that $P \cap (\Z^p \times \R^{d-p})$ is empty, or reports that $f$ is unbounded below on this set, or returns an optimal solution of $\min\{f(x) : x \in P \cap (\Z^p \times \R^{d-p})\}$.
\end{proposition}

\begin{proposition}[{\cite[Theorem~1.1]{AriHil26e}}]\label{prop:ah}
For every fixed $d$, there is a polynomial-time algorithm that, given a quadratic polynomial $f$ in $d$ variables, a polyhedron $P = \{x \in \R^d : Ax \le b\}$, and an integer $p \in \{0, \dots, d\}$, all with rational data, either reports that $P \cap (\Z^p \times \R^{d-p})$ is empty, or reports that $f$ is unbounded below on this set, or returns an optimal solution of $\min\{f(x) : x \in P \cap (\Z^p \times \R^{d-p})\}$.
\end{proposition}

The next two lemmas are standard~\cite{ConCorZam14,Ede87}. The sources do not state the polynomial running times and the bounds on the number of inequalities and of sign vectors that we need, so we include proofs in Appendix~\ref{app:prelim}. The first computes the projection of a polyhedron.

\begin{lemma}\label{lem:projection}
For every fixed $r$, there is a polynomial-time algorithm that, given a polyhedron $P = \{(x,y) \in \R^d \times \R^r : Ax + Dy \le b\}$ defined by $m$ rational inequalities, returns a rational system $Gx \le h$ with at most $(m+1)^{r+1}$ inequalities such that
\[
\{x \in \R^d : (x,y) \in P \text{ for some } y \in \R^r\} = \{x \in \R^d : Gx \le h\}.
\]
\end{lemma}

The second lists the sign vectors of an arrangement of hyperplanes in fixed dimension. Let $\psi_1, \dots, \psi_m : \R^d \to \R$ be affine functions. The \emph{sign vector} of a point $v \in \R^d$ is the vector
\[
(\sign \psi_1(v), \dots, \sign \psi_m(v)) \in \{-1, 0, 1\}^m.
\]

\begin{lemma}\label{lem:arrangement}
For every fixed $d$, there is a polynomial-time algorithm that, given affine functions $\psi_1, \dots, \psi_m : \R^d \to \R$ with rational coefficients, returns the set of sign vectors of all points of $\R^d$. This set has at most $(2m+1)^d$ elements.
\end{lemma}

\section{The optimal value of one block}\label{sec:block}

In this section we consider a single block of~\eqref{eq:bmiqp} and drop the index $i$. We treat the linking variables $x \in \R^d$ as parameters, where $d$ is arbitrary. Let $A \in \Q^{m \times d}$, $D \in \Q^{m \times r}$, $b \in \Q^m$, $R \in \Q^{d \times r}$, $c \in \Q^r$, let $Q \in \Q^{r \times r}$ be symmetric, and let $f(x,y) = y^\top Q y + 2x^\top R y + c^\top y$. The \emph{fiber} at $x \in \R^d$ is the polyhedron $P(x) = \{y \in \R^r : Dy \le b - Ax\}$. We study the optimal value $\min\{f(x,y) : y \in P(x)\}$ as a function of $x$. Proposition~\ref{prop:values} below shows that this optimal value is the minimum of polynomially many quadratic functions of $x$, each defined on a polyhedron. For fixed $r$, these functions and polyhedra can be computed in polynomial time. In Section~\ref{sec:main} we apply Proposition~\ref{prop:values} to every block of~\eqref{eq:bmiqp}.

Let $\cJ$ be the family of all sets $J \subseteq [m]$ such that the rows of $D$ indexed by $J$ are linearly independent. Every $J \in \cJ$ satisfies $|J| \le r$, so $|\cJ| \le (m+1)^r$.

\begin{proposition}\label{prop:values}
Let $r$ be fixed. For every $J \in \cJ$, one can compute in polynomial time a quadratic polynomial $g_J$ in $x$ with rational coefficients and a rational system $G_J x \le h_J$ with at most $(m + 4r + 1)^{r+1}$ inequalities, which describes a polyhedron $X_J \subseteq \R^d$, such that the following holds.
\begin{enumerate}[label=(\roman*)]
\item\label{it:val-a} For every $J \in \cJ$ and $x \in X_J$, one can compute, in time polynomial in the encoding lengths of the data and of $x$, a point $y \in P(x)$ with $f(x,y) = g_J(x)$.
\item\label{it:val-b} For every $x \in \R^d$ such that $P(x)$ is nonempty and $f(x,\cdot)$ is bounded below on $P(x)$, there is $J \in \cJ$ with $x \in X_J$, and
\[
\min\{f(x,y) : y \in P(x)\} = \min\{ g_J(x) : J \in \cJ,\ x \in X_J \}.
\]
\end{enumerate}
\end{proposition}

The fact that the optimal value is the minimum of finitely many quadratic functions, each defined on a polyhedron, is known also without a bound on $r$~\cite[Lemmas~2.4 and~2.5]{LeeTamYen05}. Proposition~\ref{prop:values} adds that, for fixed $r$, polynomially many functions suffice and that they can be computed in polynomial time.

\begin{figure}
\centering
\begin{tikzpicture}
\begin{axis}[
  width=0.62\textwidth, height=0.36\textwidth,
  axis lines=middle, xmin=-2.6, xmax=2.5, ymin=-5.5, ymax=4.8,
  xtick={-2,-1,1,2}, ytick={-4,-2,2,4},
  xlabel={$x$}, xlabel style={anchor=west},
  tick label style={font=\scriptsize},
  every axis plot/.append style={samples=60},
  clip=false]
\addplot[black!22, line width=3.2pt, domain=-2.1:0] {2*x-1};
\addplot[black!22, line width=3.2pt, domain=0:2] {x^2-1};
\addplot[black, dashed, domain=-2.1:1] {x^2};
\addplot[black, domain=-2.1:2] {x^2-1};
\addplot[black, domain=-2.1:2] {2*x-1};
\fill (axis cs:1,1) circle (1.2pt);
\fill (axis cs:2,3) circle (1.2pt);
\node[anchor=east, font=\scriptsize] at (axis cs:-2.1,4.41) {$g_{\emptyset}$};
\node[anchor=east, font=\scriptsize] at (axis cs:-2.1,3.2) {$g_{\{1\}}$};
\node[anchor=east, font=\scriptsize] at (axis cs:-2.1,-5.2) {$g_{\{2\}}$};
\end{axis}
\end{tikzpicture}
\caption{Proposition~\ref{prop:values} for the block $\min\{-y^2 + 2xy : y \ge x - 1,\ y \le 1\}$, where $d = r = 1$. The family $\cJ$ consists of $\emptyset$ and the two constraints $\{1\}$ and $\{2\}$. We have $g_{\emptyset}(x) = x^2$ on $X_{\emptyset} = \{x : x \le 1\}$, and $g_{\{1\}}(x) = x^2 - 1$ and $g_{\{2\}}(x) = 2x - 1$ on $X_{\{1\}} = X_{\{2\}} = \{x : x \le 2\}$. The optimal value of the block (shaded) is the minimum of these functions. The function $g_{\emptyset}$ comes from the stationary point $y = x$, which is a maximizer, and it is never the minimum.}
\label{fig:values}
\end{figure}
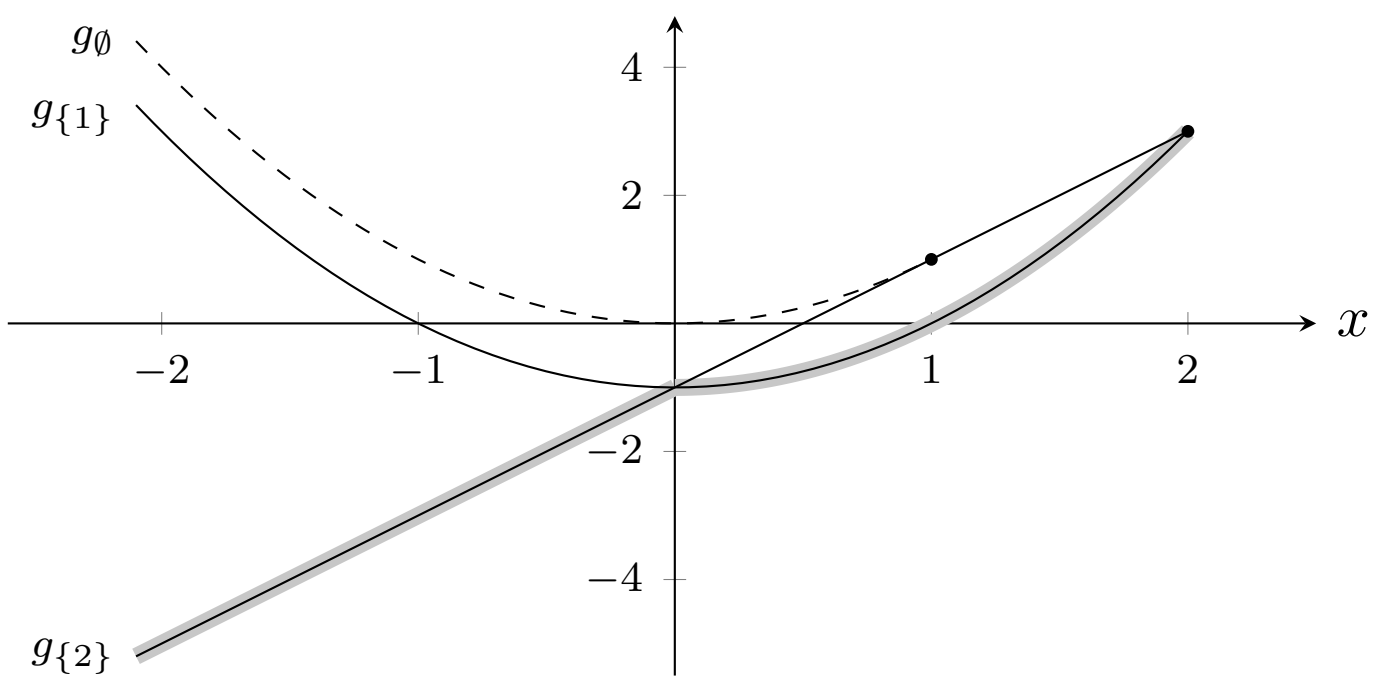

The polynomials $g_J$ and the polyhedra $X_J$ come from the stationary points of $f(x,\cdot)$ on the affine subspaces where the constraints in $J$ hold with equality. This follows the active-set argument of~\cite{Vav90} for continuous quadratic programming. That argument selects a minimizer at which the first-order conditions on the active constraints form a nonsingular linear system. Such a minimizer need not exist when the fiber contains a line, for example for the zero objective on~$\R$. We use instead that the objective is constant on the stationary points of each affine subspace, and we describe these points for all values of $x$ at once. Figure~\ref{fig:values} illustrates Proposition~\ref{prop:values} for a block with $d = r = 1$.

The gradient of $f$ with respect to $y$ is
\[
\nabla_y f(x,y) = 2Q y + 2R^\top x + c.
\]

For $J \in \cJ$ and $x \in \R^d$, let $S_J(x)$ be the set of stationary points of $f(x,\cdot)$ on the affine subspace $\{y : D_J y = b_J - A_J x\}$, that is, of the points $y \in \R^r$ with $D_J y = b_J - A_J x$ and $\nabla_y f(x,y) = -D_J^\top \lambda$ for some $\lambda \in \R^{|J|}$. Equivalently, $y \in S_J(x)$ if and only if $K_J (y, \lambda) = e_J(x)$ for some $\lambda \in \R^{|J|}$, where
\[
K_J = \begin{pmatrix} 2Q & D_J^\top \\ D_J & 0 \end{pmatrix},
\qquad
e_J(x) = \begin{pmatrix} -2R^\top x - c \\ b_J - A_J x \end{pmatrix}.
\]
The matrix $K_J$ is square of order $r + |J|$, and $e_J(x)$ is affine in $x$.

The next lemma shows that $f(x,\cdot)$ is constant on each set $S_J(x)$ and that every minimizer of $f(x,\cdot)$ over $P(x)$ lies in one of these sets.

\begin{lemma}\label{lem:stationary}
Let $x \in \R^d$.
\begin{enumerate}[label=(\roman*)]
\item\label{it:st-i} For every $J \in \cJ$, the function $f(x,\cdot)$ is constant on $S_J(x)$.
\item\label{it:st-ii} If $y^* \in P(x)$ minimizes $f(x,\cdot)$ over $P(x)$, then $y^* \in S_J(x)$ for some $J \in \cJ$.
\end{enumerate}
\end{lemma}

\begin{proof}
\ref{it:st-i}. Let $y, y' \in S_J(x)$. By the definition of $S_J(x)$, there are $\lambda, \lambda' \in \R^{|J|}$ with $\nabla_y f(x,y) = -D_J^\top \lambda$ and $\nabla_y f(x,y') = -D_J^\top \lambda'$. Let $z = y' - y$. Then $D_J z = 0$, and hence $\nabla_y f(x,y)^\top z = -\lambda^\top D_J z = 0$. Since $\nabla_y f(x,y') - \nabla_y f(x,y) = 2Qz$, subtracting the two equations for the gradients gives $2 Q z = -D_J^\top (\lambda' - \lambda)$, and hence $2 z^\top Q z = -(\lambda' - \lambda)^\top D_J z = 0$. Since $f(x,\cdot)$ is quadratic with Hessian $2Q$, we obtain $f(x,y') = f(x,y) + \nabla_y f(x,y)^\top z + z^\top Q z = f(x,y)$.

\ref{it:st-ii}. Let $I = \{j \in [m] : (D y^*)_j = (b - Ax)_j\}$, and let $J \subseteq I$ index a basis of the row space of $D_I$. Then $J \in \cJ$ and $\ker D_J = \ker D_I$. Let $z \in \ker D_J$. For every $j \in I$ and $\theta \in \R$, we have $(D(y^* + \theta z))_j = (b - Ax)_j$. For every $j \notin I$, we have $(D y^*)_j < (b - Ax)_j$. Hence $y^* + \theta z \in P(x)$ for all $\theta$ in a neighborhood of $0$. Since $y^*$ minimizes $f(x,\cdot)$ over $P(x)$, the function $\theta \mapsto f(x, y^* + \theta z)$ has a local minimum at $\theta = 0$, and its derivative there vanishes, that is, $\nabla_y f(x,y^*)^\top z = 0$. Since this holds for every $z \in \ker D_J$, the gradient $\nabla_y f(x,y^*)$ lies in the row space of $D_J$, so $\nabla_y f(x,y^*) = -D_J^\top \lambda$ for some $\lambda$. Since $J \subseteq I$, we also have $D_J y^* = b_J - A_J x$. Therefore $y^* \in S_J(x)$.
\end{proof}

The next lemma describes $S_J(x)$ by an affine function of $x$ and a subspace that does not depend on $x$. Its proof, which uses a generalized inverse of $K_J$, is in Appendix~\ref{app:prop}.

\begin{lemma}\label{lem:param}
Let $J \in \cJ$. One can compute in polynomial time rational matrices $H$, $U$, $B$ and a rational vector $\bar y$, where $B$ has $\ell \le r$ columns, such that the following holds. Let $y(x) = U x + \bar y$ and $\cL = \{B w : w \in \R^{\ell}\}$. For every $x \in \R^d$, the set $S_J(x)$ is nonempty if and only if $H e_J(x) = 0$, and in this case $S_J(x) = y(x) + \cL$.
\end{lemma}

\begin{proof}[Proof of Proposition~\ref{prop:values}]
We first define the polynomials $g_J$ and the polyhedra $X_J$. For $J \in \cJ$, let $H_J$, $U_J$, $B_J$, $\bar y_J$, $\ell_J$, and $y_J(x) = U_J x + \bar y_J$ be the objects given by Lemma~\ref{lem:param}, and let
\[
g_J(x) = f(x, y_J(x)),
\]
which is a quadratic polynomial in $x$ with rational coefficients. Let
\[
\hat X_J = \bigl\{ (x,w) \in \R^d \times \R^{\ell_J} : H_J e_J(x) = 0,\ Ax + D(y_J(x) + B_J w) \le b \bigr\},
\]
and let $X_J$ be the projection of $\hat X_J$ onto the space of the variables $x$.

Next we show that $g_J$ and $G_J x \le h_J$ can be computed as claimed. By Lemma~\ref{lem:param}, the matrices $H_J$, $U_J$, $B_J$ and the vector $\bar y_J$ can be computed in polynomial time, and hence so can $g_J$ and the system describing $\hat X_J$. The matrix $H_J$ has at most $r + |J| \le 2r$ rows. Hence the system describing $\hat X_J$ has at most $m + 4r$ inequalities, where each equation is written as two inequalities. The polyhedron $\hat X_J$ has $\ell_J \le r$ variables $w$. Hence, if $r$ is fixed, Lemma~\ref{lem:projection}, applied to $\hat X_J$ with $w$ in place of $y$, computes in polynomial time a system $G_J x \le h_J$ with at most $(m + 4r + 1)^{\ell_J + 1} \le (m + 4r + 1)^{r+1}$ inequalities that describes $X_J$.

\ref{it:val-a}. Let $J \in \cJ$ and $x \in X_J$. Since $X_J$ is the projection of $\hat X_J$, linear programming finds in polynomial time a vector $w$ with $(x,w) \in \hat X_J$. Let $y = y_J(x) + B_J w$. The inequalities in the definition of $\hat X_J$ say that $y \in P(x)$. Since $H_J e_J(x) = 0$, both $y$ and $y_J(x)$ lie in $S_J(x)$ by Lemma~\ref{lem:param}. By Lemma~\ref{lem:stationary}\ref{it:st-i}, we obtain $f(x,y) = f(x, y_J(x))$, which is $g_J(x)$ by definition.

\ref{it:val-b}. By the theorem of Frank and Wolfe~\cite[Appendix~(i)]{FraWol56}, $f(x,\cdot)$ attains its minimum $v$ on $P(x)$ at some point $y^*$. By part~\ref{it:val-a}, every $J \in \cJ$ with $x \in X_J$ satisfies $g_J(x) \ge v$. By Lemma~\ref{lem:stationary}\ref{it:st-ii}, there is $J \in \cJ$ with $y^* \in S_J(x)$. Then $S_J(x)$ is nonempty, so Lemma~\ref{lem:param} gives $H_J e_J(x) = 0$ and $y^* = y_J(x) + B_J w$ for some $w$. Since $y^* \in P(x)$, we obtain $(x,w) \in \hat X_J$, so $x \in X_J$. As in the proof of part~\ref{it:val-a}, we obtain $g_J(x) = f(x,y^*) = v$.
\end{proof}

\section{Proof of Theorem~\ref{thm:main}}\label{sec:main}

This section describes the algorithm of Theorem~\ref{thm:main} and proves its correctness and its polynomial running time. For every $i \in [n]$, let $\cJ_i$ be the family of all sets $J \subseteq [m_i]$ such that the rows of $D_i$ indexed by $J$ are linearly independent. We apply Proposition~\ref{prop:values} to block $i$, whose fiber at $x \in \R^d$ is $P_i(x) = \{y_i \in \R^{r_i} : D_i y_i \le b_i - A_i x\}$. This gives, for every $J \in \cJ_i$, a polynomial $g_{iJ}$ and a system $G_{iJ} x \le h_{iJ}$ that describes a polyhedron $X_{iJ}$. Since $r_i \le r$, we have $|\cJ_i| \le (m_i + 1)^r$, and every system $G_{iJ} x \le h_{iJ}$ has at most $(m_i + 4r + 1)^{r+1}$ inequalities.

A point $x \in \Z^p \times \R^q$ is \emph{admissible} if $(x, y_1, \dots, y_n)$ is a feasible point of~\eqref{eq:bmiqp} for some $y_1, \dots, y_n$. For every admissible $x$, every fiber $P_i(x)$ is nonempty, and its recession cone is $C_i = \{z \in \R^{r_i} : D_i z \le 0\}$. We say that block $i$ is \emph{unbounded below at} $x$ if $f_i(x,\cdot)$ is unbounded below on $P_i(x)$. If this happens for an admissible $x$, then the objective of~\eqref{eq:bmiqp} is unbounded below, since the other blocks can be fixed at feasible values. The first step of the algorithm decides whether this happens.

\begin{lemma}\label{lem:unbounded}
For all fixed $d$ and $r$, there is a polynomial-time algorithm that, given an instance of~\eqref{eq:bmiqp} that has a feasible point, decides whether some block $i$ is unbounded below at some admissible $x$.
\end{lemma}

The proof of Lemma~\ref{lem:unbounded} is in Appendix~\ref{app:unbounded}. It uses the following ray criterion. A quadratic function $f(y) = y^\top Q y + c^\top y$ is unbounded below on a nonempty polyhedron if and only if there are a feasible point $y$ and a recession direction $z$ such that either the curvature $z^\top Q z$ along $z$ is negative, or it is zero and the directional derivative $\nabla f(y)^\top z$ at $y$ along $z$ is negative. Negative curvature along a recession direction of block $i$ is detected by a quadratic program in $r_i$ variables, which can be solved in polynomial time~\cite{Vav90}, \cite[Theorem~5]{dPDeyMol17MPA}. Otherwise, the recession directions of block $i$ of zero curvature form the union of polynomially many polyhedral cones, each generated by polynomially many vectors. For every generator $z$, deciding whether the directional derivative $\nabla_y f_i(x, y_i)^\top z$ along $z$ is negative for some admissible $x$ and some $y_i \in P_i(x)$ is a mixed-integer linear program with $p$ integer variables, which Proposition~\ref{prop:milp} solves.

Once Lemma~\ref{lem:unbounded} has excluded blocks that are unbounded below at admissible points, Proposition~\ref{prop:values}\ref{it:val-b} applies to every block at every admissible $x$. The optimal value of block $i$ is then the minimum of $g_{iJ}(x)$ over the sets $J \in \cJ_i$ with $x \in X_{iJ}$. For some choices of the data, even with $d = r = 2$, the optimal value of a block is finite for every $x \in \R^d$ but is not piecewise quadratic, that is, $\R^d$ is not the union of finitely many polyhedra on each of which it is a quadratic function of $x$~\cite[Example~4.1 and Proposition~4.1]{LeeTamYen05}. Hence the set of values of $x$ at which a given function $g_{iJ}$ attains this minimum is in general not a union of finitely many polyhedra. Therefore the algorithm compares the functions $g_{iJ}$ through the signs of quadratic polynomials in $x$, which are affine functions in the space of the monomials of degree at most two.

We now state the algorithm of Theorem~\ref{thm:main}.

\begin{algo}\label{alg:main}
The input is an instance of~\eqref{eq:bmiqp}, where $d$ and $r$ are fixed.
\begin{description}[leftmargin=0pt, style=unboxed]
\item[Step~0 (feasibility and unbounded blocks).] Decide whether~\eqref{eq:bmiqp} has a feasible point by Proposition~\ref{prop:milp}, applied with zero objective to the variables $(x, y_1, \dots, y_n)$. If it does not, report that~\eqref{eq:bmiqp} has no feasible point and stop. Then apply Lemma~\ref{lem:unbounded}. If it finds a block that is unbounded below at some admissible $x$, report that the objective is unbounded below and stop. Otherwise, for every admissible $x$ and every $i \in [n]$, the function $f_i(x,\cdot)$ is bounded below on the nonempty polyhedron $P_i(x)$.

\item[Step~1 (candidates for each block).] For every $i \in [n]$ and $J \in \cJ_i$, compute $g_{iJ}$ and the system $G_{iJ} x \le h_{iJ}$ by Proposition~\ref{prop:values}.

\item[Step~2 (arrangement in the space of monomials).] Let $\bar d = d + d(d+1)/2$, and let $\nu(x) \in \R^{\bar d}$ be the vector of the monomials $x_a$ and $x_a x_b$, $1 \le a \le b \le d$. Every polynomial of degree at most two in $x$ can be written as $\alpha^\top \nu(x) + \beta$, and we associate with it the affine function $v \mapsto \alpha^\top v + \beta$ on $\R^{\bar d}$. Consider the following polynomials of degree at most two:
\begin{itemize}
\item for every inequality of every system $G_{iJ} x \le h_{iJ}$, the affine function given by its left-hand side minus its right-hand side,
\item the differences $g_{iJ} - g_{iJ'}$ for every $i \in [n]$ and $J, J' \in \cJ_i$ with $J \ne J'$.
\end{itemize}
Let $\psi_1, \dots, \psi_N : \R^{\bar d} \to \R$ be the corresponding affine functions. Lemma~\ref{lem:arrangement}, applied to $\psi_1, \dots, \psi_N$, computes the set $\Sigma$ of sign vectors of all points of $\R^{\bar d}$.

\item[Step~3 (one subproblem for each sign vector).] Let $\pi \in \Sigma$. For every $i \in [n]$, let $\cJ_i(\pi)$ be the family of sets $J \in \cJ_i$ such that $\pi$ gives a nonpositive sign to every affine function coming from $G_{iJ} x \le h_{iJ}$. If $\cJ_i(\pi)$ is empty for some $i$, discard $\pi$. Otherwise, for every $i$, choose $J_i(\pi) \in \cJ_i(\pi)$ such that $\pi$ gives a nonpositive sign to $g_{iJ_i(\pi)} - g_{iJ'}$ for every $J' \in \cJ_i(\pi)$ with $J' \neq J_i(\pi)$. Since $\pi$ is the sign vector of a point $v \in \R^{\bar d}$, such a set exists. Indeed, every set $J \in \cJ_i(\pi)$ that minimizes the value at $v$ of the affine function corresponding to $g_{iJ}$ has this property. If several sets have this property, any of them can be chosen. Consider the subproblem
\begin{equation}\label{eq:sub}
\begin{aligned}
\min\quad & x^\top Q_0 x + c_0^\top x + \sum_{i=1}^n g_{iJ_i(\pi)}(x) \\
\text{s.t.}\quad & A_0 x \le b_0, \\
& G_{iJ_i(\pi)} x \le h_{iJ_i(\pi)}, && i = 1, \dots, n, \\
& x \in \Z^p \times \R^q.
\end{aligned}
\end{equation}

\item[Step~4 (solve and recover).] Solve each subproblem~\eqref{eq:sub} by Proposition~\ref{prop:ah}. If some subproblem is unbounded below, report that the objective is unbounded below and stop. Otherwise, let $x^*$ be an optimal solution of a subproblem with smallest optimal value among the feasible subproblems, and let $\pi^*$ be the corresponding sign vector. For every $i \in [n]$, compute $y_i^* \in P_i(x^*)$ with $f_i(x^*, y_i^*) = g_{iJ_i(\pi^*)}(x^*)$ by Proposition~\ref{prop:values}\ref{it:val-a}, and return $(x^*, y_1^*, \dots, y_n^*)$.
\end{description}
\end{algo}

The correctness of Algorithm~\ref{alg:main} rests on the following relation between the subproblems and~\eqref{eq:bmiqp}.

\begin{lemma}\phantomsection\label{lem:sub}
\begin{enumerate}[label=(\roman*)]
\item\label{it:sub-a} Let $\pi \in \Sigma$ be a sign vector that Step~3 does not discard, and let $x$ be a feasible point of the subproblem~\eqref{eq:sub} of $\pi$. For every $i \in [n]$, let $y_i \in P_i(x)$ with $f_i(x, y_i) = g_{iJ_i(\pi)}(x)$, as given by Proposition~\ref{prop:values}\ref{it:val-a}. Then $(x, y_1, \dots, y_n)$ is a feasible point of~\eqref{eq:bmiqp}, and its objective value equals the objective value of $x$ in~\eqref{eq:sub}.
\item\label{it:sub-b} Assume that Step~0 does not stop. Then for every feasible point of~\eqref{eq:bmiqp}, some subproblem~\eqref{eq:sub} has a feasible point with at most the same objective value.
\end{enumerate}
\end{lemma}

\begin{proof}
\ref{it:sub-a}. Since $x$ is feasible for~\eqref{eq:sub}, we have $A_0 x \le b_0$, $x \in \Z^p \times \R^q$, and $x \in X_{iJ_i(\pi)}$ for every $i \in [n]$. Hence the vectors $y_i$ exist by Proposition~\ref{prop:values}\ref{it:val-a}. Since $y_i \in P_i(x)$ for every $i \in [n]$, the point $(x, y_1, \dots, y_n)$ is a feasible point of~\eqref{eq:bmiqp}. Its objective value equals the objective value of $x$ in~\eqref{eq:sub}, since $f_i(x, y_i) = g_{iJ_i(\pi)}(x)$ for every $i \in [n]$.

\ref{it:sub-b}. Let $(x, y_1, \dots, y_n)$ be a feasible point of~\eqref{eq:bmiqp}. Then $x$ is admissible. Since Step~0 does not stop, for every $i \in [n]$, the function $f_i(x,\cdot)$ is bounded below on the nonempty polyhedron $P_i(x)$. Let $\pi$ be the sign vector of $\nu(x)$, which belongs to $\Sigma$. Each value $\psi_j(\nu(x))$ is the value at $x$ of the corresponding polynomial, so $\pi$ gives to $\psi_j$ the sign of this polynomial at $x$. Hence, for every $i \in [n]$, a set $J \in \cJ_i$ belongs to $\cJ_i(\pi)$ if and only if $x \in X_{iJ}$. By Proposition~\ref{prop:values}\ref{it:val-b}, this family is nonempty for every $i \in [n]$, so Step~3 does not discard $\pi$. Moreover, for every $i \in [n]$,
\[
g_{iJ_i(\pi)}(x) = \min\{ g_{iJ}(x) : J \in \cJ_i(\pi) \} = \min\{f_i(x, y) : y \in P_i(x)\} \le f_i(x, y_i).
\]
The first equality holds by the choice of $J_i(\pi)$, and the second by Proposition~\ref{prop:values}\ref{it:val-b}. Since $J_i(\pi) \in \cJ_i(\pi)$, we have $x \in X_{iJ_i(\pi)}$ for every $i \in [n]$. Hence $x$ is feasible for the subproblem of $\pi$, and its value there is at most the objective value of $(x, y_1, \dots, y_n)$ in~\eqref{eq:bmiqp}.
\end{proof}

\begin{proof}[Proof of Theorem~\ref{thm:main}]
If Step~0 stops, its output is correct by Proposition~\ref{prop:milp} and Lemma~\ref{lem:unbounded}, since a block that is unbounded below at an admissible $x$ makes the objective of~\eqref{eq:bmiqp} unbounded below.

Next we assume that Step~0 does not stop, and we show that the output of Step~4 is correct. By Lemma~\ref{lem:sub}\ref{it:sub-a}, if some subproblem is unbounded below, then so is the objective of~\eqref{eq:bmiqp}. Assume now that no subproblem is unbounded below. Since~\eqref{eq:bmiqp} has a feasible point, Lemma~\ref{lem:sub}\ref{it:sub-b} shows that some subproblem is feasible. By Proposition~\ref{prop:ah}, every feasible subproblem has an optimal solution. Let $v^*$ be the smallest optimal value among the feasible subproblems. By Lemma~\ref{lem:sub}\ref{it:sub-b}, $v^*$ is at most the objective value of every feasible point of~\eqref{eq:bmiqp}. By Lemma~\ref{lem:sub}\ref{it:sub-a}, the point returned in Step~4 is a feasible point of~\eqref{eq:bmiqp} with objective value $v^*$. Hence it is an optimal solution.

Finally, we bound the running time. Since we may assume that every block has a variable, we have $n \le L$. Step~0 runs in polynomial time by Proposition~\ref{prop:milp} and Lemma~\ref{lem:unbounded}. Step~1 runs in polynomial time by Proposition~\ref{prop:values}, since there are at most $n(m+1)^r$ pairs $(i, J)$. Next we bound $N$. For every $i \in [n]$, there are at most $|\cJ_i| \le (m+1)^r$ systems $G_{iJ} x \le h_{iJ}$, each with at most $(m + 4r + 1)^{r+1}$ inequalities, and at most $|\cJ_i|^2 \le (m+1)^{2r}$ differences. Hence $N \le 2n(m + 4r + 1)^{2r+1}$, which is polynomial for fixed $r$. Since $\bar d$ is fixed, Step~2 runs in polynomial time by Lemma~\ref{lem:arrangement}, and $|\Sigma| \le (2N+1)^{\bar d} = L^{O((r+1)d^2)}$. Hence Step~3 constructs at most $|\Sigma|$ subproblems in polynomial time. Each subproblem has $d$ variables and at most $m_0 + \sum_i (m_i + 4r + 1)^{r+1}$ inequalities with polynomial encoding length, so Proposition~\ref{prop:ah} solves it in polynomial time. The optimal solution $x^*$ returned by Proposition~\ref{prop:ah} has polynomial encoding length, so Step~4 runs in polynomial time by Proposition~\ref{prop:values}\ref{it:val-a}.
\end{proof}

\section{Outline of the proof of Theorem~\ref{thm:tight}}\label{sec:tight}

We prove Theorem~\ref{thm:tight} in the following stronger form.

\begin{proposition}\label{prop:strong}
Unless P${}={}$NP, Theorem~\ref{thm:main} does not hold if any of the following changes is made, even when the linear constraints of~\eqref{eq:bmiqp} define a polytope.
\begin{enumerate}[label=(\roman*),leftmargin=*,widest=vii]
\item\label{it:s-d} The number $d$ of linking variables is not fixed, even if $q = 0$ and $n = 0$, or if $p = 0$, every block has a single variable, and the linking variables do not appear in the objective.
\item\label{it:s-k} The block size $r$ is not fixed, even if $d = 0$ and $n = 1$.
\item\label{it:s-q} The objective may contain products of variables of different blocks, even if $d = 0$ and every block has a single variable.
\item\label{it:s-link} The constraints may involve variables of different blocks, even if $d = 0$, every block has a single variable, and only two constraints involve more than one block.
\item\label{it:s-int} The blocks may contain integer variables, even if the objective is zero, $p = 3$, $q = 0$, all variables in the blocks are integral, and every block has at most $9$ variables.
\item\label{it:s-conv} There may be a single block of unbounded size whose function $f_1$ is convex in $(x, y_1)$, even if $p = 1$, $q = 0$, and the objective is $-x^2$ plus a linear function of this block.
\item\label{it:s-conv2} Blocks of unbounded size are allowed when the sum of $x^\top Q_0 x + c_0^\top x$ and of their functions $f_i$ is convex, even if $p = 1$, $q = 0$, $n = 2$, one block has a linear function $f_1$, and the other block has a single variable.
\end{enumerate}
\end{proposition}

We outline the proof of Proposition~\ref{prop:strong}, which is given in Appendix~\ref{app:tight}. In each part we reduce an NP-complete problem to deciding whether the optimal value of an instance in the stated class is at most a given number. Parts~\ref{it:s-d}--\ref{it:s-int} follow from known reductions. For part~\ref{it:s-d}, we reduce from the feasibility problem of $0/1$ integer linear programming~\cite{Kar72} and from subset sum. For parts~\ref{it:s-k} and~\ref{it:s-link}, we use the reduction of Sahni from subset sum~\cite[Theorem~2.5.4]{Sah74}. For part~\ref{it:s-q}, we reduce from maximum cut~\cite{GarJohSto76}. For part~\ref{it:s-int}, we adapt the reduction from $3$-colorability in the proof of~\cite[Theorem~21]{DvoEibGanKnoOrd21}.

Parts~\ref{it:s-conv} and~\ref{it:s-conv2} need more work. They rely on the construction of Pardalos and Vavasis~\cite{ParVav91}. For a graph $G$ and an integer $\kappa$, it gives a polyhedron in which two variables $x$ and $z$ are affine functions of the other variables. After bounding some variables, we obtain a polytope $P_2$ with the following properties. If $G$ has a clique of size $\kappa$, the minimum of $z - x^2$ over $P_2$ is $0$ and is attained at a point with $x \in \Z$. Otherwise, $P_2$ is empty or this minimum is positive. Hence declaring $x$ an integer variable keeps the minimum equal to $0$ in the clique case, and it can only increase the minimum otherwise. We take $x$ as the only linking variable and all other variables as one block. Then the objective is $-x^2$ plus the linear function $z$ of the block, which gives part~\ref{it:s-conv}. Replacing $-x^2$ by $-v^2$ for a new continuous variable $v$ with $v = x$ gives part~\ref{it:s-conv2}.

\bigskip
\begin{small}\sloppy
\paragraph{Funding.} A. Del Pia is partially supported by AFOSR Grant \mbox{FA9550-23-1-0433} and ONR Grant \mbox{N00014-25-1-2490}. Any opinions, findings, conclusions, or recommendations expressed in this material are those of the author and do not necessarily reflect the views of the Air Force Office of Scientific Research or the Office of Naval Research.

\paragraph{AI methodology.} AI tools (Claude, developed by Anthropic) were used for proofreading and suggesting improvements to the exposition. The author takes full responsibility for the content of the paper.
\end{small}

{\small
\bibliographystyle{plainurl}
\bibliography{MIQP_IPCO}
}

\clearpage
\appendix

\section{Proofs of Lemmas~\ref{lem:projection} and~\ref{lem:arrangement}}\label{app:prelim}

\begin{proof}[Proof of Lemma~\ref{lem:projection}]
Let $\Lambda = \{\lambda \in \R^m : \lambda \ge 0,\ \lambda^\top D = 0\}$ be the projection cone of $P$. Every $\lambda \in \Lambda$ gives an inequality $\lambda^\top A x \le \lambda^\top b$ that is valid for the projection of $P$ onto the space of the variables $x$, and the inequalities given by representatives of the extreme rays of $\Lambda$ describe this projection~\cite[Theorem~3.46 and Exercise~3.35]{ConCorZam14}.

Next we show that the nonzero vectors in an extreme ray of $\Lambda$ have at most $r + 1$ nonzero entries. Let $\lambda$ be a nonzero vector in an extreme ray of $\Lambda$, and let $S = \{j \in [m] : \lambda_j \neq 0\}$. For $S' \subseteq [m]$, let $\Lambda_{S'} = \{\mu \in \R^{S'} : \mu^\top D_{S'} = 0\}$. The constraints of $\Lambda$ that $\lambda$ satisfies with equality are $\lambda_j \ge 0$ for $j \notin S$ and $\lambda^\top D = 0$. Since $\Lambda$ lies in the nonnegative orthant, it is pointed, so by~\cite[Theorem~3.35]{ConCorZam14} these constraints have rank at least $m - 1$. Since $\lambda \neq 0$ satisfies them with equality, their rank is exactly $m - 1$. Hence the set of vectors that satisfy them with equality, which is isomorphic to $\Lambda_S$, is one-dimensional. Therefore $\Lambda_S$ is spanned by the restriction of $\lambda$ to $S$, whose entries are all positive. Since $\dim \Lambda_S \ge |S| - \operatorname{rank}(D_S) \ge |S| - r$, we obtain $|S| \le r + 1$.

Finally, we find a representative of every extreme ray of $\Lambda$ by enumerating the possible sets of nonzero entries. For every nonempty $S \subseteq [m]$ with $|S| \le r + 1$, we compute $\Lambda_S$. If $\Lambda_S$ is one-dimensional and spanned by a vector $\mu$ such that $\mu$ or $-\mu$ has all entries positive, we replace $\mu$ by $-\mu$ if needed and output the inequality given by the vector $\bar\mu \in \Lambda$ that equals $\mu$ on $S$ and $0$ outside $S$. Since every extreme ray of $\Lambda$ contains a vector with at most $r + 1$ nonzero entries, the output contains a representative of every extreme ray of $\Lambda$, so it describes the projection. The number of sets $S$ is at most $\sum_{j=1}^{r+1} \binom{m}{j} \le (m+1)^{r+1}$. For each $S$, Gaussian elimination computes a basis of $\Lambda_S$ and decides the remaining conditions in polynomial time, and the vector $\bar\mu$ has polynomial encoding length~\cite[Theorem~3.3]{Sch86}.
\end{proof}

\begin{proof}[Proof of Lemma~\ref{lem:arrangement}]
We first bound the number of sign vectors. We may assume that no $\psi_j$ is constant, since a constant function has the same sign at every point. The sets of points with a common sign vector are the faces of the arrangement of the distinct hyperplanes among $\{\psi_j = 0\}$, $j \in [m]$. The bound of Edelsbrunner assumes that no hyperplane is vertical, which holds after a change of coordinates. By~\cite[Theorem~1.3]{Ede87}, the number of $t$-faces of an arrangement of at most $m$ hyperplanes in $\R^d$ is at most $\sum_{a=0}^{t} \binom{d-a}{t-a} \binom{m}{d-a}$. Summing over $t$, the number of faces is at most $\sum_{j=0}^d \binom{m}{j} 2^j \le (2m+1)^d$.

Next we compute the sign vectors. We compute the sets $S_l$ of sign vectors of $\psi_1, \dots, \psi_l$, for $l = 0, 1, \dots, m$, where $S_0$ consists of the empty vector. For every $\sigma \in S_{l-1}$ and $\tau \in \{-1, 0, 1\}$, we decide whether appending $\tau$ to $\sigma$ gives a vector of $S_l$. This amounts to deciding whether a system with strict linear inequalities and linear equations has a solution. We solve, in polynomial time, the linear program that maximizes $\varepsilon$ subject to $0 \le \varepsilon \le 1$, $\psi_j(v) \le -\varepsilon$ when the sign is $-1$, $\psi_j(v) = 0$ when the sign is $0$, and $\psi_j(v) \ge \varepsilon$ when the sign is $1$. The vector belongs to $S_l$ if and only if the linear program is feasible and its optimal value is positive. By the counting bound, every set $S_l$ has at most $(2m+1)^d$ elements, so the algorithm solves polynomially many linear programs.
\end{proof}

\section{Proof of Lemma~\ref{lem:param}}\label{app:prop}

\begin{proof}[Proof of Lemma~\ref{lem:param}]
We first construct $H$ and a generalized inverse of $K = K_J$, that is, a matrix $K^-$ with $K K^- K = K$. Let $\rho$ be the rank of $K$. Let $H$ be a matrix whose rows form a basis of $\{h : h^\top K = 0\}$. Then the image of $K$ is $\{e : H e = 0\}$. Choose sets $I_1$ and $I_2$ of $\rho$ rows and $\rho$ columns of $K$ such that the submatrix $K_{I_1 I_2}$ is nonsingular. Let $K^-$ be the matrix of the same order as $K$ whose submatrix with rows indexed by $I_2$ and columns indexed by $I_1$ is $K_{I_1 I_2}^{-1}$, and whose other entries are $0$. After permuting rows and columns, we can write
\[
K = \begin{pmatrix} K_{I_1 I_2} & K_{12} \\ K_{21} & K_{22} \end{pmatrix}.
\]
Since $K$ has rank $\rho$, its columns indexed by $I_2$ span its column space. Hence there is a matrix $Y$ with $K_{12} = K_{I_1 I_2} Y$ and $K_{22} = K_{21} Y$. The first equation gives $Y = K_{I_1 I_2}^{-1} K_{12}$, and therefore $K_{22} = K_{21} K_{I_1 I_2}^{-1} K_{12}$. A direct computation then gives $K K^- K = K$. Hence, for every $e$ in the image of $K$, say $e = K u$, we have $K (K^- e) = K K^- K u = e$.

We now define $U$, $\bar y$, and $B$. Write $e_J(x) = E x + \bar e$, and let $U$ and $\bar y$ consist of the first $r$ rows of $K^- E$ and $K^- \bar e$, respectively. Let $B$ be a matrix whose columns form a basis of the subspace of the vectors formed by the first $r$ coordinates of the vectors in $\ker K$. Then $\ell \le r$.

Next we show that these matrices have the required property. If $S_J(x)$ is nonempty, then $e_J(x)$ is in the image of $K$, so $H e_J(x) = 0$. Conversely, if $H e_J(x) = 0$, then $K^- e_J(x)$ solves $K u = e_J(x)$, and the solution set of this system is $K^- e_J(x) + \ker K$. Taking the first $r$ coordinates, we obtain $S_J(x) = y(x) + \cL$.

Finally, all matrices are computed by Gaussian elimination from $K$, $E$, and $\bar e$. Therefore they are computed in polynomial time and have polynomial encoding length.
\end{proof}

\section{Proof of Lemma~\ref{lem:unbounded}}\label{app:unbounded}

We use the following three lemmas. The first is the ray criterion of Eaves.

\begin{lemma}[{\cite[Corollary~3.4]{Eav71}}]\label{lem:ray}
Let $P \subseteq \R^\ell$ be a nonempty polyhedron with recession cone $C$, and let $f(y) = y^\top Q y + c^\top y$ with $Q$ symmetric. Then $f$ is unbounded below on $P$ if and only if there are $y \in P$ and $z \in C$ such that either $z^\top Q z < 0$, or $z^\top Q z = 0$ and $\nabla f(y)^\top z < 0$.
\end{lemma}

The second lemma describes the recession directions of zero curvature.

\begin{lemma}\label{lem:zero}
Let $C = \{z \in \R^\ell : Dz \le 0\}$, and let $Q$ be a symmetric matrix with $z^\top Q z \ge 0$ for every $z \in C$. Let $\cJ$ be the family of sets $J$ of row indices of $D$ such that the corresponding rows are linearly independent. For $J \in \cJ$, let $N_J$ be a matrix whose columns form a basis of $\ker D_J$, and let
\[
Z_J = \{z \in \R^\ell : Dz \le 0,\ D_J z = 0,\ N_J^\top Q z = 0\}.
\]
Then $\{z \in C : z^\top Q z = 0\} = \bigcup_{J \in \cJ} Z_J$.
\end{lemma}

\begin{proof}
We first show that $Z_J \subseteq \{z \in C : z^\top Q z = 0\}$ for every $J \in \cJ$. Let $z \in Z_J$. Since $N_J^\top Q z = 0$, the vector $Qz$ is orthogonal to $\ker D_J$, so $Qz = D_J^\top \mu$ for some $\mu$. Hence $z^\top Q z = \mu^\top D_J z = 0$. Moreover, $z \in C$, since $Dz \le 0$.

Next we show the reverse inclusion. Let $z^* \in C$ with $(z^*)^\top Q z^* = 0$. Since $z^\top Q z \ge 0$ for every $z \in C$, the vector $z^*$ minimizes $z^\top Q z$ over $C$. Let $I$ be the set of indices $j$ with $(D z^*)_j = 0$, and let $J \subseteq I$ index a basis of the row space of $D_I$. Then $J \in \cJ$, $D_J z^* = 0$, and $\ker D_J = \ker D_I$. Since $(D z^*)_j < 0$ for every $j \notin I$, for every $z \in \ker D_J$ the vector $z^* + \theta z$ lies in $C$ for all $\theta$ in a neighborhood of $0$. Hence the derivative of $\theta \mapsto (z^* + \theta z)^\top Q (z^* + \theta z)$ at $\theta = 0$ vanishes, that is, $z^\top Q z^* = 0$. Therefore $N_J^\top Q z^* = 0$ and $z^* \in Z_J$.
\end{proof}

The third lemma computes a finite set of generators of a polyhedral cone. It is standard~\cite{ConCorZam14,Sch86}. The sources do not state the polynomial running time and the bound on $|\Gamma|$ that we need, so we include a proof.

\begin{lemma}\label{lem:generators}
For every fixed $\ell$, there is a polynomial-time algorithm that, given a rational matrix $A$ with $m$ rows and $\ell$ columns, returns a finite set $\Gamma$ of rational vectors in $C = \{z \in \R^\ell : Az \le 0\}$, with $|\Gamma| \le 2\ell + (m+1)^{\ell-1}$, such that every vector of $C$ is a nonnegative combination of vectors of $\Gamma$.
\end{lemma}

\begin{proof}
We first construct $\Gamma$. Let $B$ be a matrix whose columns form a basis of the lineality space $\ker A$ of $C$, and put every column of $B$ and its negative into $\Gamma$. The cone $C' = \{z \in C : B^\top z = 0\}$ is pointed, and $C = \ker A + C'$. Hence it suffices to add to $\Gamma$ a representative of every extreme ray of $C'$~\cite[Section~8.8]{Sch86}. For every set $S$ of at most $\ell - 1$ rows of $A$, compute the kernel of the matrix obtained by stacking $A_S$ and $B^\top$. If this kernel is one-dimensional and spanned by a vector $\rho$ such that $\rho$ or $-\rho$ lies in $C'$, add that vector to $\Gamma$.

Next we show that $\Gamma$ contains a representative of every extreme ray of $C'$. Let $z$ be a nonzero vector in an extreme ray of $C'$. By~\cite[Theorem~3.35]{ConCorZam14}, the constraints of $C'$ that $z$ satisfies with equality have rank $\ell - 1$. These constraints include $B^\top z = 0$. Hence there is a set $S$ of at most $\ell - 1$ rows of $A$ with $A_S z = 0$ such that the matrix obtained by stacking $A_S$ and $B^\top$ has rank $\ell - 1$. The kernel of this matrix is spanned by $z$, so the algorithm adds a positive multiple of $z$ to $\Gamma$.

Finally, we bound $|\Gamma|$ and the running time. The matrix $B$ has at most $\ell$ columns, which give at most $2\ell$ vectors. There are at most $(m+1)^{\ell-1}$ sets $S$, and each adds at most one vector, so $|\Gamma| \le 2\ell + (m+1)^{\ell-1}$. Each kernel is computed in polynomial time by Gaussian elimination.
\end{proof}

\begin{proof}[Proof of Lemma~\ref{lem:unbounded}]
We first describe the pairs $(x, y_i)$ with $x$ admissible and $y_i \in P_i(x)$ by linear inequalities and integrality constraints. For every $j \in [n]$, let $X_j = \{x \in \R^d : P_j(x) \neq \emptyset\}$. Lemma~\ref{lem:projection}, applied to the polyhedron $\{(x, y_j) \in \R^d \times \R^{r_j} : A_j x + D_j y_j \le b_j\}$, computes a system of at most $(m_j + 1)^{r+1}$ linear inequalities that describes $X_j$, since $r_j \le r$. Fix $i \in [n]$, and let $F_i$ be the set of pairs $(x, y_i)$ such that $x \in \Z^p \times \R^q$, $A_0 x \le b_0$, $x \in X_j$ for every $j \ne i$, and $y_i \in P_i(x)$. Then $(x, y_i) \in F_i$ if and only if $x$ is admissible and $y_i \in P_i(x)$. By Lemma~\ref{lem:ray}, applied to $P = P_i(x)$ and $f = f_i(x,\cdot)$ for every admissible $x$, block $i$ is unbounded below at some admissible $x$ if and only if there are $(x, y_i) \in F_i$ and $z \in C_i$ such that either $z^\top Q_i z < 0$, or $z^\top Q_i z = 0$ and $\nabla_y f_i(x, y_i)^\top z < 0$. The set $F_i$ is nonempty, since~\eqref{eq:bmiqp} has a feasible point.

Next we decide whether $z^\top Q_i z < 0$ for some $z \in C_i$. Solve the quadratic program $\min\{z^\top Q_i z : z \in C_i,\ -1 \le z \le 1\}$ in $r_i \le r$ variables by Proposition~\ref{prop:ah} with $p = 0$ (here the polynomial-time solvability of quadratic programming in fixed dimension suffices~\cite{Vav90}, \cite[Theorem~5]{dPDeyMol17MPA}). By homogeneity, its optimal value is negative if and only if $z^\top Q_i z < 0$ for some $z \in C_i$. In this case, block $i$ is unbounded below at every admissible $x$ by Lemma~\ref{lem:ray}.

Next we treat the case where $z^\top Q_i z \ge 0$ for every $z \in C_i$. Then Lemma~\ref{lem:zero}, applied with $D = D_i$ and $Q = Q_i$, shows that the directions $z \in C_i$ with $z^\top Q_i z = 0$ form the union of the cones $Z_{iJ}$ over $J \in \cJ_i$. For every such $J$, $Z_{iJ}$ is described by $m_i + 2r_i$ homogeneous inequalities, since $D_J$ has $|J|$ rows, $N_J$ has $r_i - |J|$ columns, and each equation is written as two inequalities. Lemma~\ref{lem:generators}, applied to this system with $\ell = r_i$, computes a set $\Gamma_{iJ}$ of generators of $Z_{iJ}$. There are $|\cJ_i| \le (m_i + 1)^r$ such sets $J$. Every vector of $Z_{iJ}$ is a nonnegative combination of vectors of $\Gamma_{iJ} \subseteq Z_{iJ}$. Since $\nabla_y f_i(x, y_i)^\top z$ is linear in $z$, it is negative for some $z \in Z_{iJ}$ if and only if it is negative for some $\rho \in \Gamma_{iJ}$. Hence block $i$ is unbounded below at some admissible $x$ if and only if, for some $J \in \cJ_i$ and some $\rho \in \Gamma_{iJ}$, the infimum of the linear function
\[
(x, y_i) \mapsto \nabla_y f_i(x, y_i)^\top \rho = (2 Q_i y_i + 2 R_i^\top x + c_i)^\top \rho
\]
over $F_i$ is negative or $-\infty$. Each of these problems is a mixed-integer linear program with $p$ integer variables, which Proposition~\ref{prop:milp} solves in polynomial time. This decides whether block $i$ is unbounded below at some admissible $x$.

Finally, we bound the running time. The number of blocks, sets $J$, and generators is polynomial, and every quadratic program and mixed-integer linear program above has polynomial encoding length. Hence the algorithm runs in polynomial time.
\end{proof}

\section{Proof of Proposition~\ref{prop:strong}}\label{app:tight}

Part~\ref{it:s-int} relies on the reduction from $3$-colorability in the proof of~\cite[Theorem~21]{DvoEibGanKnoOrd21}. We write it with inequalities instead of slack variables, which gives smaller blocks. Since Proposition~\ref{prop:strong} concerns polytopes, we also add bounds that make the feasible region a polytope. The next lemma states the properties that we need.

\begin{lemma}\label{lem:col}
Given a graph $G$, one can construct in polynomial time a system of linear constraints in variables $c_1, c_2, c_3$ and further variables with the following properties.
\begin{enumerate}[label=(\roman*),leftmargin=*]
\item\label{it:col-a} The system has an integral solution if and only if $G$ has a $3$-coloring.
\item\label{it:col-b} The feasible region of the system is a polytope.
\item\label{it:col-c} The further variables can be partitioned into groups of at most $9$ variables such that every constraint involves only $c_1, c_2, c_3$ and the variables of at most one group.
\end{enumerate}
\end{lemma}

\begin{proof}
We first construct the system. Let $[t]$ be the vertex set of $G$, and let $\mathrm{pr}(i)$ be the $i$th prime. For a variable $c$, an integer $s \ge 2$, and variables $m$, $r$, $u$, let $\mathrm{Div}(c, s; m, r, u)$ denote the constraints
\begin{equation}\label{eq:col}
c = s\, m + r, \quad m \ge 0, \quad 0 \le r \le s - 1, \quad u \le r \le (s - 1)\, u, \quad 0 \le u \le 1.
\end{equation}
In every integral solution of~\eqref{eq:col}, $r$ is the remainder of $c$ modulo $s$, since $0 \le r \le s - 1$. Moreover, $u \in \{0,1\}$, so the constraints $u \le r \le (s - 1)\, u$ give $u = 0$ if and only if $r = 0$. Hence $u = 0$ if and only if $s$ divides $c$. The system consists of the constraints
\begin{equation}\label{eq:colsys}
\begin{aligned}
& \mathrm{Div}(c_j, \mathrm{pr}(i); m_{ij}, r_{ij}, u_{ij}), && i \in [t],\ j \in [3],\\
& u_{i1} + u_{i2} + u_{i3} = 2, && i \in [t],\\
& \mathrm{Div}(c_j, \mathrm{pr}(i); m_{eij}, r_{eij}, u_{eij}), && e = \{i,k\} \in E(G),\ j \in [3],\\
& \mathrm{Div}(c_j, \mathrm{pr}(k); m_{ekj}, r_{ekj}, u_{ekj}), && e = \{i,k\} \in E(G),\ j \in [3],\\
& u_{eij} + u_{ekj} \ge 1, && e = \{i,k\} \in E(G),\ j \in [3],\\
& 0 \le c_j \le \mathrm{pr}(1) \cdots \mathrm{pr}(t), && j \in [3].
\end{aligned}
\end{equation}
The first two lines ensure that $\mathrm{pr}(i)$ divides exactly one of $c_1, c_2, c_3$, for every $i$. The third, fourth, and fifth lines ensure that $\mathrm{pr}(i)$ and $\mathrm{pr}(k)$ do not both divide $c_j$, for every edge $\{i,k\}$ and every $j$. Each edge has its own copies of the variables of $i$ and $k$, because otherwise the constraint $u_{eij} + u_{ekj} \ge 1$ would join the variables of $i$ and $k$ into one group in~\ref{it:col-c}. Without the constraints $c_j \le \mathrm{pr}(1) \cdots \mathrm{pr}(t)$, this is the system in the proof of~\cite[Theorem~21]{DvoEibGanKnoOrd21}, where every equation with a nonnegative slack variable is written as an inequality. By the prime number theorem, $\mathrm{pr}(t) = O(t \log t)$. Hence the product $\mathrm{pr}(1) \cdots \mathrm{pr}(t) \le \mathrm{pr}(t)^t$ has encoding length $O(t \log t)$, and the system can be constructed in polynomial time.

Next we prove~\ref{it:col-a}. Dvo\v{r}\'ak, Eiben, Ganian, Knop, and Ordyniak showed that their system has an integral solution if and only if $G$ has a $3$-coloring~\cite[proof of Theorem~21]{DvoEibGanKnoOrd21}. Eliminating the slack variables does not change this, since the slack variables are integral whenever the other variables are. The constraints $c_j \le \mathrm{pr}(1) \cdots \mathrm{pr}(t)$ can only remove solutions. They keep the solution that Dvo\v{r}\'ak et al.\ construct from a $3$-coloring, since in that solution $c_j$ is the product of the primes $\mathrm{pr}(i)$ over the vertices $i$ of color $j$, which is at most $\mathrm{pr}(1) \cdots \mathrm{pr}(t)$.

Next we prove~\ref{it:col-b} and~\ref{it:col-c}. The variables $c_j$ are bounded by the last line of~\eqref{eq:colsys}, the variables $r$ and $u$ are bounded by~\eqref{eq:col}, and every $m$ is bounded because $s\, m = c - r$ in~\eqref{eq:col}. Hence the feasible region is a polytope. For~\ref{it:col-c}, we take one group for every vertex $i$, consisting of the $9$ variables $m_{ij}, r_{ij}, u_{ij}$ with $j \in [3]$, and one group for every edge $e$ and $j \in [3]$, consisting of the $6$ variables indexed by $e$ and $j$. The constraints in the first two lines of~\eqref{eq:colsys} for $i$ involve only $c_1, c_2, c_3$ and the group of $i$. The constraints in the third, fourth, and fifth lines for $e$ and $j$ involve only $c_j$ and the group of $e$ and $j$. The last line involves no group.
\end{proof}

Parts~\ref{it:s-conv} and~\ref{it:s-conv2} rely on a construction of Pardalos and Vavasis~\cite{ParVav91}, which reduces the clique problem to the minimization of $z - x^2$ over a polyhedron. This polyhedron may be unbounded. Since Proposition~\ref{prop:strong} concerns polytopes, we add bounds that make it a polytope. The proof of part~\ref{it:s-conv} also needs that, in the clique case, the minimum is attained at a point with $x \in \Z$. The next lemma states the properties that we need.

\begin{lemma}\label{lem:pv}
Given a graph $G$ and a positive integer $\kappa$, one can construct in polynomial time a polytope $P_2$ of points $(\xi, \eta, x, z)$, where $x, z \in \R$, with the following properties.
\begin{enumerate}[label=(\roman*),leftmargin=*]
\item\label{it:pv-a} If $G$ has a clique of size $\kappa$, then the minimum of $z - x^2$ over $P_2$ is $0$ and is attained at a point with $x \in \Z$.
\item\label{it:pv-b} If $G$ has no clique of size $\kappa$, then $P_2$ is empty or the minimum of $z - x^2$ over $P_2$ is positive.
\end{enumerate}
\end{lemma}

\begin{proof}
We first construct $P_2$. Let $[t]$ be the vertex set of $G$. Consider the variables $\xi_1, \dots, \xi_t$, $\eta_{jl}$ for $1 \le j < l \le t$, $x$, and $z$, and the constraints
\begin{equation}\label{eq:pv}
\begin{aligned}
& 0 \le \xi_j \le 1, && j = 1, \dots, t, \\
& x = \sum_{j=1}^t 4^j \xi_j, \\
& \eta_{jl} \ge 0,\quad \eta_{jl} \ge \xi_j + \xi_l - 1, && 1 \le j < l \le t, \\
& z = \sum_{j=1}^t 4^{2j} \xi_j + 2 \sum_{1 \le j < l \le t} 4^{j+l} \eta_{jl}, \\
& \xi_j + \xi_l \le 1, && \{j,l\} \notin E(G), \\
& \xi_1 + \dots + \xi_t = \kappa .
\end{aligned}
\end{equation}
These are the constraints of~\cite[Sections~2 and~3]{ParVav91} with the constant $b = 4$, where we write $\xi_j$, $\eta_{jl}$, and $x$ for their variables $x_j$, $y_{jl}$, and $w$. Let $P_0$ be the polyhedron defined by the constraints in~\eqref{eq:pv} except the last two lines, let $P_1$ be the polyhedron defined by all constraints in~\eqref{eq:pv}, and let $P_2$ be the polyhedron defined by~\eqref{eq:pv} together with $\eta_{jl} \le 1$ for $1 \le j < l \le t$. All variables are bounded on $P_2$, since $0 \le \xi_j \le 1$ and $0 \le \eta_{jl} \le 1$, and $x$ and $z$ are determined by $\xi$ and $\eta$. Hence $P_2$ is a polytope. The coefficients $4^j$ and $4^{j+l}$ have encoding length $O(t)$, so $P_2$ can be constructed in polynomial time. Pardalos and Vavasis proved that the minimum of $z - x^2$ over $P_0$ is $0$~\cite[Lemma~1 and Theorem~1]{ParVav91}. They also proved that the minimum of $z - x^2$ over $P_1$ is $0$ if and only if $G$ has a clique of size $\kappa$~\cite[Section~3]{ParVav91}. Since $P_2 \subseteq P_1 \subseteq P_0$, the function $z - x^2$ is nonnegative on $P_1$.

Next we prove~\ref{it:pv-a}. Let $K$ be a clique of size $\kappa$. Let $\xi_j = 1$ for $j \in K$ and $\xi_j = 0$ otherwise, let $\eta_{jl} = \xi_j \xi_l$, and define $x$ and $z$ by~\eqref{eq:pv}. Pardalos and Vavasis showed that this point belongs to $P_1$ and that $z - x^2 = 0$ at this point~\cite[Lemma~1 and Section~3]{ParVav91}. Since $\eta_{jl} \le 1$, the point belongs to $P_2$. Moreover, $x = \sum_{j \in K} 4^j$ is an integer. Hence the minimum of $z - x^2$ over $P_2$ is $0$ and is attained at this point.

Next we prove~\ref{it:pv-b}. Assume that $G$ has no clique of size $\kappa$ and that $P_2$ is nonempty. Since $P_2$ is a polytope, $z - x^2$ attains its minimum over $P_2$. If this minimum were $0$, then the minimum of $z - x^2$ over $P_1$ would be $0$, since $z - x^2 \ge 0$ on $P_1$ and $P_2 \subseteq P_1$. By the result of Pardalos and Vavasis, $G$ would then have a clique of size $\kappa$, a contradiction. Hence the minimum of $z - x^2$ over $P_2$ is positive.
\end{proof}

\begin{proof}[Proof of Proposition~\ref{prop:strong}]
In each part we give a polynomial-time reduction from an NP-complete problem to deciding whether the optimal value of an instance in the stated class is at most a given number. All feasible regions are polytopes.

\ref{it:s-d}. We first prove the statement with $q = 0$ and $n = 0$. Deciding whether a system $A_0 x \le b_0$, $0 \le x \le 1$, has an integral solution is NP-complete~\cite{Kar72}. This is the feasibility question for~\eqref{eq:bmiqp} with zero objective, $p$ integer variables, $q = 0$, and no blocks.

Next we prove the statement with $p = 0$. Let $a_1, \dots, a_t$ and $\beta$ be positive integers. Consider the instance with continuous linking variables $x_1, \dots, x_t$, the linking constraints $\sum_{j=1}^t a_j x_j = \beta$ and $0 \le x_j \le 1$ for $j \in [t]$, and, for every $j \in [t]$, a block with a single variable $y_j$, the constraint $y_j = x_j$, and the function $f_j(x, y_j) = y_j - y_j^2$. Here $p = 0$, every block has a single variable, the linking variables do not appear in the objective, and the feasible region is a polytope. Since $y_j = x_j$, the objective equals $\sum_{j=1}^t (x_j - x_j^2)$, which is nonnegative on $[0,1]^t$ and vanishes exactly at the points of $\{0,1\}^t$. Hence the optimal value is $0$ if and only if some subset of $\{a_1, \dots, a_t\}$ sums to $\beta$, and this subset sum problem is NP-complete~\cite{Kar72}.

\ref{it:s-k}. Let $a_1, \dots, a_t$ and $\beta$ be positive integers, and consider
\[
\min\Bigl\{ \sum_{j=1}^t (y_j - y_j^2) : \sum_{j=1}^t a_j y_j = \beta,\ 0 \le y_j \le 1,\ j = 1, \dots, t \Bigr\},
\]
with all variables in one block. Here $d = 0$ and $n = 1$, and the size of the block is not bounded. Every term $y_j - y_j^2$ is nonnegative on $[0,1]$ and vanishes exactly at $0$ and $1$. Hence the optimal value is $0$ if and only if the equation has a $0/1$ solution, that is, if and only if some subset of $\{a_1, \dots, a_t\}$ sums to $\beta$. This subset sum problem is NP-complete~\cite{Kar72}. This is essentially the reduction of Sahni~\cite[Theorem~2.5.4]{Sah74}.

\ref{it:s-q}. Let $G$ be a graph with vertex set $[n]$, and consider
\[
\min\Bigl\{ -\sum_{\{i,j\} \in E(G)} (y_i + y_j - 2 y_i y_j) : 0 \le y_i \le 1,\ i = 1, \dots, n \Bigr\}.
\]
Here $d = 0$, every variable forms its own block, and the objective contains products of variables of different blocks. The objective is affine in each variable separately, so it attains its minimum over the box at a vertex. At a vertex $y \in \{0,1\}^n$, the term $y_i + y_j - 2y_iy_j$ equals $1$ if the edge $\{i,j\}$ crosses the cut defined by $y$ and $0$ otherwise. Hence, for every integer $k$, the optimal value is at most $-k$ if and only if $G$ has a cut with at least $k$ edges, and deciding this is NP-complete~\cite{GarJohSto76}.

\ref{it:s-link}. Consider the problem of part~\ref{it:s-k}, but now let block $j$ consist of the single variable $y_j$, for every $j \in [t]$. Here $d = 0$, $n = t$, every block has a single variable, and the equation gives the only two constraints that involve more than one block. The optimal value is the same as in part~\ref{it:s-k}.

\ref{it:s-int}. Let $G$ be a graph, and consider the system of Lemma~\ref{lem:col} with zero objective. We take $c_1, c_2, c_3$ as integer linking variables, so $p = 3$ and $q = 0$. By Lemma~\ref{lem:col}\ref{it:col-b}, the feasible region is a polytope. By Lemma~\ref{lem:col}\ref{it:col-c}, the groups of the other variables form integral blocks with at most $9$ variables. By Lemma~\ref{lem:col}\ref{it:col-a}, the instance is feasible if and only if $G$ has a $3$-coloring. Since $3$-colorability is NP-complete~\cite{GarJohSto76}, this proves part~\ref{it:s-int}.

\ref{it:s-conv}. Let $G$ be a graph, let $\kappa$ be a positive integer, and let $P_2$ be the polytope of Lemma~\ref{lem:pv}. Consider the problem $\min\{z - x^2 : (\xi,\eta,x,z) \in P_2,\ x \in \Z\}$. Here $x$ is the only linking variable, so $p = 1$ and $q = 0$. The variables $\xi$, $\eta$, $z$ form one block. The objective is $-x^2$ plus the linear function $z$ of the block, so $f_1$ is convex in $(x, y_1)$. The optimal value is at least the minimum of $z - x^2$ over $P_2$. By Lemma~\ref{lem:pv}, if $G$ has a clique of size $\kappa$, the optimal value is $0$, and otherwise the problem is infeasible or its optimal value is positive. Since the clique problem is NP-complete~\cite{Kar72}, this proves part~\ref{it:s-conv}.

\ref{it:s-conv2}. Modify the instance of part~\ref{it:s-conv} by adding a continuous variable $v$, the constraint $v = x$, and replacing the objective by $z - v^2$. The linking variable $x$ is integral, the variables $\xi$, $\eta$, $z$ form a block with the linear function $f_1 = z$, and $v$ forms a block of size one. The linking part of the objective is zero, so its sum with $f_1$ is linear, hence convex. Since $v = x$ on the feasible region, the optimal value is the same as in part~\ref{it:s-conv}.

\end{proof}

\section{Proof of Proposition~\ref{thm:xp}}\label{app:xp}

\begin{proof}[Proof of Proposition~\ref{thm:xp}]
We first consider the instances with $n = 0$, parameterized by $d$. If $n = 0$ and $q = 0$, problem~\eqref{eq:bmiqp} is IQP in $d$ variables. Herrmann proved that deciding whether the optimal value of an IQP is at most a given number is W[1]-hard parameterized by the number of variables~\cite{Her26}. The hardness persists for instances over polytopes. Indeed, Del Pia, Dey, and Molinaro showed that if an MIQP has a feasible point with objective value at most a given number, then it has one whose encoding length is bounded by a polynomial in the input~\cite{dPDeyMol17MPA}. Adding the corresponding bounds on all variables gives an equivalent instance over a polytope with the same number of variables, in polynomial time. An algorithm with running time $\omega(d,r) \cdot L^c$ would decide these instances in fixed-parameter time, and hence FPT${}={}$W[1].

Next we consider the instances with $d = 0$ and $n = 1$, parameterized by $r$. If $d = 0$ and $n = 1$, problem~\eqref{eq:bmiqp} is a continuous quadratic program in $r$ variables. Deciding whether $z^\top Q z < 0$ for some $z$ with $Az \le 0$ is W[1]-hard parameterized by the number of variables~\cite[Theorem~6.6]{AriHil26c}. By homogeneity, this holds if and only if the optimal value of $\min\{z^\top Q z : Az \le 0,\ -1 \le z \le 1\}$ is negative, and this problem is an instance with $d = 0$, $n = 1$, and $r$ equal to the number of variables, over a polytope. An algorithm with running time $\omega(d,r) \cdot L^c$ would decide these instances in fixed-parameter time, and hence FPT${}={}$W[1].
\end{proof}

\section{Proof of Corollary~\ref{cor:vi}}\label{app:cor}

\begin{proof}[Proof of Corollary~\ref{cor:vi}]
We first find the blocks. Since the interaction graph has vertex integrity at most $v$, there is a set $W$ of continuous variables such that $|W|$ plus the largest number of nodes in a connected component after deleting $W$ is at most $v$. We find such a set by enumerating all sets of at most $v$ continuous variables and computing the connected components after deleting each of them. This takes polynomial time for fixed $v$.

Next we write the MIQP in the form~\eqref{eq:bmiqp}. Declare the integer variables and the variables in $W$ linking, and take the connected components as blocks. Two continuous variables in different components are not adjacent. Hence every constraint and every term of the objective involves the linking variables and at most one block. The constraints that involve only linking variables form $A_0 x \le b_0$, and each other constraint belongs to the system of its block. The terms of the objective split in the same way into $Q_0$, $c_0$, $R_i$, $Q_i$, and $c_i$. We obtain an instance of~\eqref{eq:bmiqp} with $d \le p + v$ and $r \le v$, which Theorem~\ref{thm:main} solves in polynomial time.
\end{proof}

\section{Theorem~\ref{thm:main} from integer quadratic programming}\label{app:wei}

Proposition~\ref{prop:ah} is the only result of~\cite{AriHil26e} used in the proof of Theorem~\ref{thm:main}. Here we show that it can be replaced by the following result of Wei, which concerns integer quadratic programs.

\begin{proposition}[{\cite[Theorem~1.1]{Wei26}}]\label{prop:wei}
For every fixed $d$, there is a polynomial-time algorithm that, given a quadratic polynomial $f$ in $d$ variables and a polyhedron $P = \{x \in \R^d : Ax \le b\}$, all with rational data, either reports that $P \cap \Z^d$ is empty, or reports that $f$ is unbounded below on this set, or returns an optimal solution of $\min\{f(x) : x \in P \cap \Z^d\}$.
\end{proposition}

We first show that Theorem~\ref{thm:main} holds for instances with $q = 0$ when Proposition~\ref{prop:ah} is replaced by Proposition~\ref{prop:wei}. Proposition~\ref{prop:ah} is used in two places. The proof of Lemma~\ref{lem:unbounded} uses it only with $p = 0$, where it is the polynomial-time solvability of quadratic programming in fixed dimension~\cite{Vav90}, \cite[Theorem~5]{dPDeyMol17MPA}. Step~4 of Algorithm~\ref{alg:main} and the proof of Theorem~\ref{thm:main} use it to solve the subproblems~\eqref{eq:sub}. If $q = 0$, every subproblem~\eqref{eq:sub} is an integer quadratic program in $p$ variables, which Proposition~\ref{prop:wei} solves, with the same three outcomes as Proposition~\ref{prop:ah}.

Next we derive Proposition~\ref{prop:ah} from this special case. Consider an MIQP in $d$ variables, $p$ of them integral. If $p = 0$, it is a quadratic program in fixed dimension, so assume $p \ge 1$. Let $x$ consist of the integer variables and $y$ of the continuous ones. Up to an additive constant, the objective has the form $x^\top Q_0 x + c_0^\top x + y^\top Q_1 y + 2x^\top R_1 y + c_1^\top y$, and the constraints have the form $A_1 x + D_1 y \le b_1$. This is an instance of~\eqref{eq:bmiqp} with $q = 0$, linking variables $x$, and one block $y$ of size $d - p$. Since $p \le d$ and $d - p \le d$ are fixed, Theorem~\ref{thm:main} for $q = 0$ solves it in polynomial time, with the outcomes required by Proposition~\ref{prop:ah}.

Hence Theorem~\ref{thm:main} also follows from Proposition~\ref{prop:wei} instead of Proposition~\ref{prop:ah}. In particular, the polynomial-time solvability of MIQP in fixed dimension~\cite{AriHil26e} also follows from that of integer quadratic programming in fixed dimension~\cite{Wei26}.

\end{document}